\documentclass[a4paper,12pt]{article}
\usepackage{amssymb}

\newcommand{\R}{\mathbb{R}}

\newcommand{\beq}{\begin{equation} }
\newcommand{\eqq}{\end{equation} }
\newcommand{\cuad}{{\sqcap\kern-.68em\sqcup}}

\newcommand{\equ}[1]{(\ref{#1})}

\newtheorem{teo}{Theorem}[section]

\newtheorem{proposition}{Proposition}[section]

\newtheorem{lemma}{Lemma}[section]

\newtheorem{remark}{Remark}[section]
\newcommand{\bremark}{\begin{remark} \em}
\newcommand{\eremark}{\end{remark} }

\def\beeq{\begin{equation}}
\def\eeq{\end{equation}}
\newcommand{\begeqaet}{\begin{eqnarray*}}
\newcommand{\eneqaet}{\end{eqnarray*}}

\begin{document}
\begin{center}{\bf    A Liouville type theorem for Lane-Emden systems involving the fractional Laplacian}\medskip

\bigskip

\bigskip

{Alexander Quaas and Aliang Xia}

Departamento de  Matem\'atica,  Universidad T\'ecnica Federico Santa Mar\'{i}a

Casilla: V-110, Avda. Espa\~na 1680, Valpara\'{\i}so, Chile.

 {\sl (alexander.quaas@usm.cl and aliangxia@gmail.com)}
\end{center}

\bigskip

\begin{abstract}
 We establish a Liouville type theorem for the fractional
Lane-Emden system:
 \begin{eqnarray*}
\left\{\begin{array}{l@{\quad }l}
(-\Delta)^\alpha u=v^q&{\rm in}\,\,\R^N,\\
(-\Delta)^\alpha v=u^p&{\rm in}\,\,\R^N,
 \end{array}
 \right.
 \end{eqnarray*}
 where $ \alpha\in(0,1) $, $ N>2\alpha $ and $ p,q $ are positive real numbers and in an appropriate new range.
 To prove our result we will use the local realization of fractional Laplacian, which can be constructed as Dirichlet-to-Neumann operator of a degenerate elliptic equation in the spirit of Caffarelli and Silvestre \cite{CS}. Our proof is based on a monotonicity argument for suitable transformed
 functions and the method of moving planes in an infinity half cylinder based on some maximum principles which obtained by some barrier functions and a coupling argument using fractional Sobolev trace inequality.
\end{abstract}
\date{}

\setcounter{equation}{0}
\section{ Introduction}

During the last years there has been a renewed and increasing interest in
the study of linear and nonlinear integral operators, including the fractional
Laplacian, motivated by many applications and by important advances on
the theory of nonlinear partial differential equations.
This because these types of systems appear as limiting equations of many phenomena, such as pattern formation, population evolution, chemical reaction, etc.
Some of these equations are named by Lotka-Volterra, Bose-Einstein, Schr\"odinger system, Gierer-Meinhardt. The solutions in most of case represent concentrations in the process and thus naturally positive solutions of the systems are of particular interest.

Most of the results in this field are obtained assuming that the diffusion is governed by the Laplacian or other more general local elliptic operators.
The mathematical literature in the case of elliptic systems when the diffusion is governed by L\'evy stable process and the elliptic operator turns to be the fractional order of Laplacian  is very recent. See \cite{DP}, \cite{MC}, \cite{QX}, \cite{ST1}, \cite{ST2},  \cite{ST3} and \cite{Wei}. Notice that the fractional Laplacian appears in different diffusion models, see for example  \cite{A}, \cite{Be}, \cite{BoG}, \cite{H}, \cite{NPV}, \cite{Ta}  and the references therein.

As far as we know, there aren't existence results for systems without variation structure. When the variational structure breaks, the methods developed to prove existence results for local elliptic systems  are obtained by Perron's Method or topological arguments; for example, the Leray-Schauder degree or Krasnoselskii's index theory, where the many assumptions in using these theories are the  priori bound for  solutions. These bounds are obtained in many settings by the classical scaling or blow-up argument due to Gidas and Spruck \cite{GS} in the scalar case and \cite{DS} for systems case, see also  references therein.
Liouville type Theorems are crucial to get a contradiction for the limiting system or equation.  Roughly speaking,  better Liouville type Theorems give more general existence results.
Observe that   there are still some problems even in the case ($\alpha=1$) which is known as Lane-Emden conjecture, see  \cite{SZ2}, \cite{S1} and blow.

The aim of this paper is establish a new Liouville type theorem for the following Lane-Emden system involving the fractional Laplacian:

\begin{eqnarray}\label{1.1}
\left\{\begin{array}{l@{\quad }l}
(-\Delta)^\alpha u=v^q&{\rm in}\,\,\R^N,\\
(-\Delta)^\alpha v=u^p&{\rm in}\,\,\R^N.
 \end{array}
 \right.
 \end{eqnarray}
where $\alpha\in(0,1)$ and $ N>2\alpha $.

The fractional Laplacian $ (-\Delta)^\alpha $ can be defined, for example, by the Fourier transform. Namely, for a function $ u $ in the Schwartz class $ \mathcal{S} $, we have
\[\widehat{(-\Delta)^\alpha} u(\xi)=|\xi|^{2\alpha}\hat{u}(\xi).\]
Furthermore, consider the space
\[
\mathcal{L}_\alpha(\R^N):=\left\{u\,\,\bigg\vert\,\,u:\R^N\rightarrow\R\,\,{\rm sucht\,\,that}\,\, \int_{\R^N}\frac{|u(y)|}{1+|y|^{N+2\alpha}}dy<\infty \right\}.
\]
endowed with the norm
\[
\|u\|_{\mathcal{L}_\alpha(\R^N)}:=\int_{\R^N}\frac{|u(y)|}{1+|y|^{N+2\alpha}}dy<\infty.
\]
If $ u\in \mathcal{L}_\alpha(\R^N) $ (see \cite{SL}), then $ (-\Delta)^\alpha u$ can be defined as a distribution, that is, for any $ \varphi\in \mathcal{S} $,
\[
\int_{\R^N} \varphi(-\Delta)^\alpha u dx=\int_{\R^N} u(-\Delta)^\alpha \varphi dx.
\]
In additional, for some $ \sigma>0 $, suppose that $ u\in \mathcal{L}_\alpha(\R^N)\cap C^{2\alpha+\sigma} (\R^N)$ if $ 0<\alpha<1/2 $ or $ u\in \mathcal{L}_\alpha(\R^N)\cap C^{1,2\alpha+\sigma-1} (\R^N)$ if $ \alpha\ge1/2 $, then we have
\[
(-\Delta)^\alpha u(x)=C_{N,\alpha} P.V.\int_{\R^N}\frac{u(x)-u(y)}{|x-y|^{N+2\alpha}}dy \quad {\rm for}\,\,x\in \R^N,
\]
where $P.V.$ denotes the principal value of the integral and $ C_{N,\alpha} $ is a normalization constant.

\medskip

When $ \alpha=1 $, the Lane-Emden system for Laplace operator
\begin{eqnarray}\label{1.5}
\left\{\begin{array}{l@{\quad }l}
-\Delta u=v^q&{\rm in}\,\,\R^N,\\
-\Delta v=u^p&{\rm in}\,\,\R^N,
 \end{array}
 \right.
 \end{eqnarray}
has been extensively studied in the literature, see \cite{BM,FF,M,SZ1,SZ2}. It has been conjectured that the Sobolev's hyperbola
\[\left\{p>0,q>0:\frac{1}{p+1}+\frac{1}{q+1}=1-\frac{2}{N}\right\},\]
is the dividing curve between existence and nonexistence for (\ref{1.5}).
For the radial case, this conjecture was completely solved by \cite{M,SZ1}. In fact, if
the pair $ (p,q) $ lies below  Sobolev's hyperbola, that is,
\beq\label{1.6}
\frac{1}{p+1}+\frac{1}{q+1}>1-\frac{2}{N},
\eqq
then there are no radial positive solutions to system (\ref{1.5}), see \cite{M} (for $ p>1,q>1 $) and \cite{SZ1} (for $ p>0,q>0 $). In addition, there are indeed positive radial solutions to system (\ref{1.5}) if $ (p,q) $ lies above  Sobolev's hyperbola (see also \cite{SZ1}).

The conjecture for more general case, i.e., without the assumption of radial symmetry, has not been completely answered yet. Partial results for nonexistence are known.
Define
\beq\label{1.9}
\gamma_1=\frac{2(q+1)}{pq-1},\quad \gamma_2=\frac{2(p+1)}{pq-1}\quad ({\rm if}\,\,pq>1).
\eqq
There are no positive classical supersolutions to (\ref{1.5}) if
\beq\label{1.10}
pq\le1\quad {\rm or}\quad pq>1 \,\,{\rm and}\,\,\max\{\gamma_1,\gamma_2\}\ge N-2,
\eqq
see \cite{M,S,SZ2}. Moreover, we know that the condition (\ref{1.9}) is  optimal for supersolution.
For positive solution,  Felmer and Figueiredo \cite{FF} proved that if
\beq\label{1.11}
0<p,q\le \frac{N+2}{N-2},\quad (p,q)\neq\left(\frac{N+2}{N-2},\frac{N+2}{N-2}\right),
\eqq
then problem (\ref{1.5}) has no classical positive solutions. Notice that for $ N\ge3 $, condition (\ref{1.9}) and (\ref{1.11}) are stronger than (\ref{1.6}).
The full conjecture seem difficult for nonradial solutions. As far as we know the full conjecture is true when $ N=3,4 $, see \cite{SZ2} and \cite{S1}. In the high dimensions, apart from (\ref{1.10}),  the conjecture is only known to be true in  some subregion of subcritial range:
\beq\label{1.8}
\min\{\gamma_1,\gamma_2\}\ge\frac{N-2}{2}\quad {\rm and}\quad(\gamma_1,\gamma_2)\neq\left(\frac{N-2}{2}, \frac{N-2}{2} \right)
\eqq
by Busca and Man\'{a}sevich \cite{BM}. Note that the condition (\ref{1.11}) in particular contains where both exponents are subcritical, that  is, the region considered in \cite{FF}.

The aim of the present paper is to show that the result of Busca and Man\'{a}sevich \cite{BM}
can be extended to system (\ref{1.1}). We prove the following result.

\begin{teo}\label{t}
Let $p,q>0$ and $pq>1$ and set
\beq\label{beta}
\beta_1=\frac{2\alpha(q+1)}{pq-1},\quad\beta_2=\frac{2\alpha(p+1)}{pq-1}.
\eqq
Suppose
\beq\label{2.3}
\beta_1,\beta_2\in\left[ \frac{N-2\alpha}{2},\,\, N-2\alpha\right)\quad and\quad (\beta_1,\beta_2)\not=\left( \frac{N-2\alpha}{2},\,\, \frac{N-2\alpha}{2}\right).
\eqq
Then, for some $ \sigma>0 $, there exists no positive $ \mathcal{L}_\alpha(\R^N)\cap C^{2\alpha+\sigma} (\R^N)$ if $ 0<\alpha<1/2 $ or  in $\mathcal{L}_\alpha(\R^N)\cap C^{1,2\alpha+\sigma-1} (\R^N)$ if $ \alpha\ge1/2 $ type solution to system (\ref{1.1}).
\end{teo}

\begin{remark}
1) Observe  that region (\ref{2.3}) in particular contains where both exponents are  subcritical, that is
\beq\label{pq}
\frac{N}{N-2\alpha}<p,q\le\frac{N+2\alpha}{N-2\alpha},\quad  with\quad(p,q)\neq\left(\frac{N+2\alpha}{N-2\alpha},\frac{N+2\alpha}{N-2\alpha}\right).
\eqq
By Theorem 3 in \cite{MC} and Theorem 3 in \cite{ZCCY}, we know there are no  positive   $ u,v\in\mathcal{L}_\alpha(\R^N)\cap C^{2\alpha+\sigma} (\R^N)$ if $ 0<\alpha<1/2 $ or  in $ \mathcal{L}_\alpha(\R^N)\cap C^{1,2\alpha+\sigma-1} (\R^N)$ if $ \alpha\ge1/2 $ type solution to system (\ref{1.1}) if (\ref{pq}) holds. Hence,  Theorem  \ref{t} is valid for a large region of $ (p,q) $ in comparison with  Theorem 3 in \cite{MC}.\\

\noindent 2) If we take $q=1$ we can obtained Liouville type resutlts for bi-fractional equations.

\end{remark}

\begin{remark}
We note that region (\ref{2.3}) does not include the point
\[
(\beta_1,\beta_2)=\left(\frac{N-2\alpha}{2}, \frac{N-2\alpha}{2}\right).
\]
Indeed, if $ \beta_1=\beta_2=\frac{N-2\alpha}{2} $, then
\[
p=q=\frac{N+2\alpha}{N-2\alpha},
\]
and problem
\beq\label{ee}
(-\Delta)^\alpha u=u^{\frac{N+2\alpha}{N-2\alpha}}.
\eqq
has nontrivial nonnegative solutions called fractional bubble, see Chen-Li-Ou \cite{CLO,CLO1}, Jin-Li-Xiong \cite{JLX} and also Y.Y. Li \cite{Li}.
\end{remark}

 In \cite{CS}, Caffarelli and Silvestre   introduced a local realization of the fractional Laplacian $ (-\Delta)^\alpha $ in $ \R^N $ through the Dirichlet-Neumann map of an appropriate degenerate elliptic operator in upper half space $ \R^{N+1}_+ $. More precisely, consider an extension of  $ u $ to the upper half space $ \R^{N+1}_+ $ so that $ U(x,0)=u(x) $ and
\[
\Delta_x U +\frac{1-2\alpha}{y}U_y+ U_{yy}=0\quad{\rm for}\,\,X=(x,y)\in \R^{N+1}_+.
\]
Let $ P_\alpha(x,y) $ denote the corresponding Poisson kernel
\[
P_\alpha(x,y)=c_{N,\alpha}\frac{y^{2\alpha}}{(|x|^2+y^2)^{(N+2\alpha)/2}}\quad {\rm for}\,\,x\in\R^N\,\,{\rm and}\,\,y>0.
\]
where $ c_{N,\alpha} $ is a normalization constant ( for an explicit value of $ c_{N,\alpha} $  see \cite{CaS}).
If $ u\in \mathcal{L}_\alpha(\R^N) $, we can define
\[
U(x,y)=P_\alpha(\cdot,y)\ast u=c_{N,\alpha}\int_{\R^N}\frac{y^{2\alpha}}{(|x-\xi|^2+y^2)^{(N+2\alpha)/2}}u(\xi)d\xi.
\]
Moreover, for some $ \sigma>0 $, suppose that $ u\in \mathcal{L}_\alpha(\R^N)\cap C^{2\alpha+\sigma} (\R^N)$ if $ 0<\alpha<1/2 $ or $ u\in \mathcal{L}_\alpha(\R^N)\cap C^{1,2\alpha+\sigma-1} (\R^N)$ if $ \alpha\ge1/2 $, then $ U\in C^2(\R^{N+1}_+)\cap C(\overline{\R^{N+1}_+}) $, $ y^{1-2\alpha}\partial_y U\in C(\overline{\R^{N+1}_+}) $ and
\begin{eqnarray*}
\left\{\begin{array}{l@{\quad }l}
div(y^{1-2\alpha}\nabla U)=0&{\rm in}\,\,\R^{N+1}_+,\\
 U=u&{\rm on}\,\,\partial\R^{N+1}_+,\\
 -\lim_{y\rightarrow0^+}y^{1-2\alpha}U_y=\kappa_\alpha (-\Delta)^\alpha u&{\rm on}\,\,\partial\R^{N+1}_+,
 \end{array}
 \right.
 \end{eqnarray*}
where
\[
\kappa_\alpha=\frac{\Gamma(1-\alpha)}{2^{2\alpha-1}\Gamma(\alpha)}
\]
with $ \Gamma $ being the Gamma function,  see Theorem 1.3 in \cite{DDW} and also \cite{CS,CaS,FF,S}.

Using the local formulation established by Caffarelli and Silvestre  \cite{CS}, the above theorem will follow as a corollary of the following  Liouville type result for a degenerated systems with a coupling with a nonlinear Neumann condition in the upper half space $ \R^{N+1}_+ $.

\begin{teo}\label{t1}
Let $p,q>0,pq>1 $ and (\ref{2.3}) holds. Then there exists no positive $ C^2(\R^{N+1}_+)\cap C(\overline{\R^{N+1}_+}) $  and $  y^{1-2\alpha}\partial_y (\cdot)\in C(\overline{\R^{N+1}_+}) $ type solution of
\begin{eqnarray}\label{ex}
\left\{\begin{array}{l@{\quad }l}
div(y^{1-2\alpha}\nabla U)=0&{\rm in}\,\,\R^{N+1}_+,\\
 -\lim_{y\rightarrow0^+}y^{1-2\alpha}U_y=V^q&{\rm on}\,\,\partial\R^{N+1}_+,\\
 div(y^{1-2\alpha}\nabla V)=0&{\rm in}\,\,\R^{N+1}_+,\\
 -\lim_{y\rightarrow0^+}y^{1-2\alpha}V_y=U^p&{\rm on}\,\,\partial\R^{N+1}_+.
 \end{array}
 \right.
 \end{eqnarray}
\end{teo}

Our proof follows the idea in \cite{BM}. Roughly speaking, as in \cite{BM}, we first transform the elliptic equation (\ref{ex}) in upper half space $ \R^{N+1}_+ $ to an  appropriate equation in upper half infinite cylinder $ \R\times S^N_+ $ (see (\ref{2.6}) and (\ref{2.9})), where $ S^N_+ $ is the upper half unit sphere. Then, we study the nonexistence result via a symmetry and monotonicity result (i.e., Lemma \ref{l2}) which obtained by the method of moving planes. However, there are some
difficulties appear compared our article with \cite{BM} that the operator  in \equ{ex} is a degenerated operator and the nonlinearity is at the boundary.
In particular, to prove the maximum principle for "narrow" domains which permits us to get the moving planes started. For this purpose we follow some similar arguments as in
\cite{FW} which are for the single equation and a coupling argument establishd by a fractional Sobolev trace inequality  (see Lemma \ref{3.1}).
We also need prove two Hopf's lemmas where barrier functions need to be construct, see Lemmas \ref{h1} and \ref{p1} in section 2.

\medskip
We end the introduction by mention that we can use Theorem  \ref{t} to obtain a priori estimate and existence result
for positive solutions of nonlinear  elliptic systems involving the fractional Laplacian.

\medskip
The paper is organized as follows. In section 2, we do a transformation as in \cite{BM} to problem (\ref{ex}) and present some preliminary results, the Hopf's lemmas and the strong maximum principle. A monotonicity  and symmetry result is shown by the method of moving planes in section 3 and we prove the nonexistence result (Theorem \ref{t1}) at the end of section 3.

\setcounter{equation}{0}
\section{ Preliminaries}

This section is devote to introduce some preliminary results, the Hopf's lemmas and the strong maximum principle. We start this section by transforming (\ref{ex}) as in  \cite{BM} by using polar coordinates and Emden-Fowler variables.
We take standard polar coordinates in $ \R^{N+1}_+: $ $ X=(x,y)=r\theta $, where $ r=|X| $ and $ \theta=X/|X| $. Denote $ \theta=(\theta_1,\theta_2,\cdots,\theta_N,\theta_{N+1}) $ and let $ \theta_{N+1}=y/|X|  $ denote the component of $ \theta $
in the $ y $ direction and $ S^N_+=\{X\in \R^{N+1}_+: r=1, \theta_{N+1}>0\} $ denote the  upper unit half space.

For a given function $ w  $ of $ X\in\R^{N+1}_+ $, we write, using the same symbol $ w $ without risk of confusion,
\[w(X)=w(r,\theta).\]
Thus we have the following formula
\beq\label{f}
\Delta w=\frac{1}{r^2} \Delta_\theta w+\frac{N}{r}\frac{\partial w}{\partial r}+\frac{\partial^2 w}{\partial r^2},
\eqq
where $ \Delta_\theta $ denotes the Laplace-Beltrami operator on $ S^N $.

Set
\begin{eqnarray}\label{2.5}
\left\{\begin{array}{l@{\quad }l}
\overline{U}(t,\theta)=r^{\beta_1}U(r,\theta)\\
\overline{V}(t,\theta)=r^{\beta_2}V(r,\theta)
 \end{array}
 \right.
 \end{eqnarray}
for $ \beta_1,\beta_2 $ to be fixed, where $ t=\log r $. Using the formula (\ref{f}), an easy calculation verifies that
\begin{eqnarray}\label{2.6}
\left\{\begin{array}{l@{\quad }l}
\theta_{N+1}^{2\alpha-1}div(\theta_{N+1}^{1-2\alpha}\nabla_\theta \overline{U})+\overline{U}_{tt}-\delta_1\overline{U}_t
-\nu_1\overline{U}=0&{\rm in}\,\, \R\times S_+^N,\\
 -\lim_{\theta_{N+1}\rightarrow 0^+}\theta_{N+1}^{1-2\alpha}\partial_{\theta_{N+1}} \overline{U}=r^{\beta_1+2\alpha-q\beta_2}\overline{V}^q&{\rm on}\,\,\R\times\partial S^N_+,\\
\theta_{N+1}^{2\alpha-1}div(\theta_{N+1}^{1-2\alpha}\nabla_\theta \overline{V})+\overline{V}_{tt}-\delta_2\overline{V}_t
-\nu_2\overline{V}=0&{\rm in}\,\, \R\times S^N_+,\\
 -\lim_{\theta_{N+1}\rightarrow 0^+}\theta_{N+1}^{1-2\alpha}\partial_{\theta_{N+1}} \overline{V}=r^{\beta_2+2\alpha-p\beta_1}\overline{U}^p&{\rm on}\,\, \R\times\partial S^N_+,
 \end{array}
 \right.
 \end{eqnarray}
where
\begin{eqnarray}\label{2.7}
\left\{\begin{array}{l@{\quad }l}
\delta_1=2\beta_1-(N-2\alpha),\quad \nu_1=\beta_1((N-2\alpha)-\beta_1),\\
\delta_2=2\beta_2-(N-2\alpha),\quad \nu_2=\beta_2((N-2\alpha)-\beta_2).
 \end{array}
 \right.
 \end{eqnarray}
For ease of the notation, we define the operators
\begin{eqnarray*}
L_\alpha \overline{U}:=\theta_{N+1}^{2\alpha-1}div(\theta_{N+1}^{1-2\alpha}\nabla \overline{U})&=&\Delta_\theta \overline{U}+\frac{1-2\alpha}{\theta_{N+1}}\frac{\partial \overline{U}}{\partial \theta_{N+1}}\\
L_\alpha \overline{V}:=\theta_{N+1}^{2\alpha-1}div(\theta_{N+1}^{1-2\alpha}\nabla \overline{V})&=&\Delta_\theta \overline{V}+\frac{1-2\alpha}{\theta_{N+1}}\frac{\partial \overline{U}}{\partial \theta_{N+1}}.
\end{eqnarray*}
If we define $ \beta_1 $ and $ \beta_2 $ as in (\ref{beta}), then we write (\ref{2.6}) as
\begin{eqnarray}\label{2.9}
\left\{\begin{array}{l@{\quad }l}
L_\alpha \overline{U}+\overline{U}_{tt}-\delta_1\overline{U}_t
-\nu_1\overline{U}=0&{\rm in}\,\, \R\times S^N_+,\\
 -\lim_{\theta_{N+1}\rightarrow 0^+}\theta_{N+1}^{1-2\alpha}\partial_{\theta_{N+1}} \overline{U}=\overline{V}^q&{\rm on}\,\, \R\times \partial S^N_+,\\
L_\alpha \overline{V}+\overline{V}_{tt}-\delta_2\overline{V}_t
-\nu_2\overline{V}=0&{\rm in}\,\, \R\times S^N_+,\\
 -\lim_{\theta_{N+1}\rightarrow 0^+}\theta_{N+1}^{1-2\alpha}\partial_{\theta_{N+1}} \overline{V}=\overline{U}^p&{\rm on}\,\, \R\times \partial S^N_+.
 \end{array}
 \right.
 \end{eqnarray}
Here we have used the facts $ \beta_1+2\alpha-q\beta_2=0 $ and $ \beta_2+2\alpha-p\beta_1 =0$ by (\ref{beta}). Moreover, with these notations in (\ref{2.7}), the assumptions (\ref{2.3})  is equivalent to
\begin{eqnarray}\label{2.8}
\left\{\begin{array}{l@{\quad }l}
\delta_1, \delta_2\ge 0,\quad (\delta_1, \delta_2)\not=(0,0),\\
\nu_1,\nu_2>0,\\
p,q>0,\quad pq>1.
 \end{array}
 \right.
 \end{eqnarray}

 In order to prove Theorem \ref{t1}, we will use the method of moving planes. The key tools for  use the method of moving planes is the Hopf's lemma and the strong maximum principle. The remain of the section is devote to prove these results related the operators we studied. We first show the following weak maximum principle.

\begin{proposition}\label{pw}
Let $ \Omega $ be an bounded domain in $ \R\times S^N_+ $ and $ w\in C^2(\Omega)\cap C(\bar{\Omega}) $. Suppose
\[
L_\alpha w+w_{tt}+a(t,\theta) w_t\le0\quad in \,\,\Omega,
\]
where $ |a(t,\theta)|\le a_0$=constant in $\Omega$.
 Then the nonnegative minimum of $ w $ in  $ \bar{\Omega} $ is achieved on $ \partial \Omega $, that is,
\[
\inf_{\Omega}w=\inf_{\partial\Omega}w.
\]
\end{proposition}

{\bf Proof.} It is clear that if $ L_\alpha w+w_{tt}+a(t,\theta)w_t<0 $ in $ \Omega $, then a strong maximum principle holds, that is, $ w $ cannot achieve an interior nonnegative minimum in $ \bar{\Omega} $. Indeed, if $ (t_0,\theta_0)\in \Omega $, then
\[
w_{tt}(t_0,\theta_0)\ge0,\quad \Delta_\theta w(t_0,\theta_0)\ge0,\quad  w_t(t_0,\theta_0)=0,\quad{\rm and}\quad\nabla_\theta w(t_0,\theta_0)=0.
\]
This implies $ L_\alpha w(t_0,\theta_0)+w_{tt}(t_0,\theta_0)+a(t_0,\theta_0)w_t(t_0,\theta_0)\ge0 $ which is impossible.

A simple computation, for $\gamma>0$,  gives,
\[
 (e^{\gamma t})_{tt}+a(e^{\gamma t})_{t}=e^{\gamma t}(\gamma^2 +a\gamma)\ge e^{\gamma t}(\gamma^2 -a_0\gamma).
\]
So we can choose $ \gamma $ large enough such that $ (e^{\gamma t})_{tt}+a(e^{\gamma t})_{t}>0 $.
Hence, for any $ \varepsilon>0 $,
\[
L_\alpha (w-\varepsilon e^{\gamma t}) +(w-\varepsilon e^{\gamma t})_{tt}+a (w-\varepsilon e^{\gamma t})_t<0
\]
in $ \Omega $ so that
\[
\inf_{\Omega} (w-\varepsilon e^{\gamma t})=\inf_{\partial \Omega}(w-\varepsilon e^{\gamma t}).
\]
 Letting $ \varepsilon\rightarrow0 $, we see that
\[
\inf_{\Omega}w=\inf_{\partial\Omega}w.
\]
as asserted in the proposition. $\Box$\\

Next we suppose more generally that
\[
L_\alpha w+w_{tt}+a(t,\theta)w_t-b(t,\theta) w\le0\quad{\rm in }\,\,\Omega,
\]
where $ |a(t,\theta)|\le a_0$ and $b$ is a nonnegative function in $\Omega$.
 By considering  the sunset $ \Omega^-\subset \Omega $ in which $ w<0 $. We can observe that
if $ L_\alpha w+w_{tt}+a(t,\theta) w_t-b(t,\theta) w\le0 $ in $ \Omega $, then  $ L_\alpha w+w_{tt}+a(t,\theta)w_t\le b(t,\theta)w\le0 $ in $ \Omega^- $ and hence the minimum of $ w $ on $ \overline{\Omega^-} $ must be achieved on $ \partial\overline{\Omega^-} $ and hence also on $ \partial\Omega $. Thus, writing $ w^-=\min\{w,0\} $ we obtain:

\begin{teo}\label{tw}
Let $ \Omega $ be an bounded domain in $ \R\times S^N_+ $ and $ w\in C^2(\Omega)\cap C(\bar{\Omega}) $. Suppose
\[
L_\alpha w+w_{tt}+a(t,\theta) w_t-b(t,\theta) w\le0\quad in \,\,\Omega,
\]
where $ |a(t,\theta)|\le a_0$ and $b$ is a nonnegative function in $\Omega$. Then
\[
\inf_{\Omega}w\ge\inf_{\partial\Omega}w^-.
\]
\end{teo}

Next, we prove two Hopf's Lemmas.  Let $ \mathcal{D} $ is an bounded domain of $ S^N_+ $ and $ \bar{t}\in \R $, we define
\[
\Omega_\delta=\{|t-\bar{t}|\le\delta\}\times \mathcal{D}\subset \R\times S^N_+
\]
for some $ \delta>0 $. The first Hopf's lemma is

\begin{lemma}\label{h1}
Suppose that $ (t_0,\theta_0)\in \partial \Omega_\delta $ and
 $ w\in C^2(\Omega_\delta)\cap C(\Omega_\delta\cup (t_0,\theta_0)) $ be a solution of
\[
L_\alpha w+w_{tt}+a(t,\theta) w_t-b(t,\theta)w\le0\quad in\,\,\Omega_\delta,
\]
where $a$ and $b$ are bounded functions and $b$ is nonnegative. Assume in addition that $ w(t,\theta)>0 $ for every  $ (t,\theta)\in\overline{\Omega}_\delta$ and $t\not=t_0$. Moreover, $ w(t_0,\theta)=0 $ if $ (t_0,\theta)\in\overline{\Omega}_\delta$. Then
\[
\lim_{t\rightarrow t_0}\frac{w(t,\theta)-w(t_0,\theta_0)}{t-t_0}<0.
\]
\end{lemma}

{\bf Proof.} For $ 0<\rho<\delta$,
 we define an auxiliary function $ \phi $  as
\[
\phi(t)=e^{-\beta |t-\bar{t}|^2}-e^{-\beta \delta^2},
\]
where $|t-\bar{t}|>\rho$  and $ \beta $ is a positive constant to be determined later. We notice that $ 0<\phi(t)<1 $.

A direct calculation gives
\begin{eqnarray*}
\phi_{tt}+a \phi_t-b \phi =e^{-\beta |t-\bar{t}|^2}\left(4\beta^2 |t-\bar{t}|^2+2\beta (a|t-\bar{t}|- 1)-b\right).
\end{eqnarray*}
Hence we can chose $ \beta $ large enough such that
\[
\phi_{tt}+a \phi_t-b \phi\ge0
\]
in $\Omega_\delta\setminus  \Omega_{\delta/2}$. Since $ w>0 $ in $ (t,\theta)\in\overline{\Omega}_\delta$ and $t\not=t_0$ and $ \phi(t_0)=0 $,  there is a $ \varepsilon>0 $ such that  $ w-\varepsilon\phi\ge0 $ on $ \partial \Omega_\delta\cup \partial \Omega_{\delta/2} $.  Moreover, we have
\[
L_\alpha (w-\varepsilon\phi)+(w-\varepsilon\phi)_{tt}+c_1 (w-\varepsilon\phi)_t-c_2 (w-\varepsilon\phi)\le 0
\]
in $\Omega_\delta\setminus  \Omega_{\delta/2}$. Hence, by the weak maximum principle (see Theorem \ref{tw}) implies that
\[
w-\varepsilon\phi\ge0 \quad{\rm in}\,\, \Omega_\delta\setminus  \Omega_{\delta/2}.
\]
Taking the outer normal derivative at $ (t_0,\theta_0) $, we obtain
\[
\lim_{t\rightarrow t_0}\frac{w(t,\theta)-w(t_0,\theta_0)}{t-t_0}\le
 \varepsilon\lim_{t\rightarrow t_0}\frac{\phi(t)-\phi(t_0)}{t-t_0}=-2\varepsilon\beta \delta e^{-\beta \delta^2}<0,
\]
as we required. $ \Box $\\

\begin{remark}
If in addition $w\in C^1(\Omega_\delta\cup (t_0,\theta_0))$, then we have
\[
\frac{\partial w}{\partial t}(t_0,\theta_0)<0.
\]
\end{remark}

Next we established the second Hopf's lemma  on the boundary $ \R\times \partial S^N_+ $.
We denote as before that $ \theta=(\tilde{\theta}, \theta_{N+1})\in S^N_+ $ with $ \theta_{N+1}>0 $.  For $ (\tilde{\theta}_0,0)\in \partial S^N_+ $, we define
\[
\Gamma_R^0=\{\tilde{\theta}\in S^{N-1}\,\,|\,\, d_{S^{N-1}}(\tilde{\theta},\tilde{\theta}_0)\le R\},
\]
as a neighbourhood of $ (\tilde{\theta}_0,0) $ on $ \partial S^N_+ $,   where $ d_{S^{N-1}} $ denotes the distance in $ S^{N-1} $. Then, for some cosntant $\tau>0$, we let
\[
\mathcal{C}_{R,\tau}(t_0,\tilde{\theta}_0,0):=(t_0-\tau, t_0+\tau)\times \Gamma_R^0\subset  \R\times \partial S^N_+.
\]
For notational simplicity we denote $ \mathcal{C}_{R,\tau}=\mathcal{C}_{R,\tau}(t_0,\tilde{\theta}_0,0) $ in what follows.

\begin{lemma}\label{p1}
Assume $ (t_0,\tilde{\theta}_0,0)\in\R\times \partial S^N_+ $ and
consider the subset $ \mathcal{C}_{R,\tau}\times (0,\varrho)$ of $ \R\times S^{N}_+ $ for $ \tau >0$ and $ 0<\varrho<1 $. Let $ w\in C^2(\mathcal{C}_{R,\tau}\times (0,\varrho))\cap C(\overline{\mathcal{C}_{R,\tau}\times (0,\varrho)}) $ satisfies
\begin{eqnarray*}
\left\{\begin{array}{l@{\quad }l}
L_\alpha w+w_{tt}+a(t,\theta) w_t-b(t,\theta) w\le0 &{ in}\,\,\overline{\mathcal{C}_{R,\tau}\times (0,\varrho)},\\
w>0&{ in}\,\,\overline{\mathcal{C}_{R,\tau}\times (0,\varrho)},\\
w(t_0,\tilde{\theta}_0,0)=0,
 \end{array}
 \right.
 \end{eqnarray*}
where $a$ and $b$ are bounded nonnegative functions.

 Then,
 \[
-\limsup_{\theta_{N+1}\rightarrow 0^+}\theta_{N+1}^{1-2\alpha}\frac{w(t_0,\tilde{\theta}_0,\theta_{N+1})}{\theta_{N+1}}<0.
 \]
\end{lemma}

{\bf Proof.} Here we follow the argument in \cite{CaS}.
Consider the function $\phi= \phi(t)$  on $ \mathcal{C}_{R,\tau} $ and satisfies the following ODE
\begin{eqnarray*}
\left\{\begin{array}{l@{\quad }l}
-\phi_{tt}-\underline{a} \phi_t +\overline{b}\phi =0 &{ \rm in}\,\,(t_0-\tau/2, t_0+\tau/2),\\
\phi>0 &{ \rm in}\,\,(t_0-\tau/2, t_0+\tau/2),\\
\phi(t_0-\tau/2)=0,
 \end{array}
 \right.
 \end{eqnarray*}
 where $\underline{a}=\inf_{\mathcal{C}_{R,\tau}\times (0,\varrho)} a\ge0$ and $\overline{b}=\sup_{\mathcal{C}_{R,\tau}\times (0,\varrho)} b\ge0$. Since $\overline{b}\ge0$,
 then we can write
 \[
 \phi(t)=C_1e^{\left(-\underline{a}+\sqrt{\underline{a}^2+4\overline{b}}\right)t}+C_2e^{\left(-\underline{a}-\sqrt{\underline{a}^2+4\overline{b}}\right)t}
 \]
with $C_1\ge0 \ge C_2$. This implies $\phi_t\ge0$ in $(t_0-\tau/2, t_0+\tau/2)$. Moreover, we can choose $ \|\phi\|_{L^\infty}\le C $.

Hence, we have
\begin{eqnarray*}
\left\{\begin{array}{l@{\quad }l}
L_\alpha \phi+\phi_{tt}+a \phi_t-b \phi=(a-\underline{a})\phi_t+(\overline{b}-b)\phi\ge0 &{ \rm in}\,\,\mathcal{C}_{R/2,\tau/2}\times (0,\varrho),\\
\phi\ge 0&{\rm in}\,\,\mathcal{C}_{R/2,\tau/2}\times (0,\varrho),\\
\phi=0,&{\rm on}\,\,\partial\mathcal{C}_{R/2,\tau/2}\times [0,\varrho),
 \end{array}
 \right.
 \end{eqnarray*}

Therefore, for $ \varepsilon>0 $,
\[
L_\alpha (w-\varepsilon \phi)+(w-\varepsilon \phi)_{tt}+a (w-\varepsilon \phi)_t-b (w-\varepsilon \phi)\le0 \quad { \rm in}\,\,\mathcal{C}_{R/2,\tau/2}\times (0,\varrho)
\]
and $ w-\varepsilon \phi=w\ge0 $ on $ \partial\mathcal{C}_{R/2,\tau/2}\times [0,\varrho) $. Moreover, taking $ \varepsilon>0 $ small enough, we have $w\ge \varepsilon \phi  $ on $ \mathcal{C}_{R/2,\tau/2}\times (\{\theta_{N+1}=\varrho/2\} \cup \{\theta_{N+1}=0\})$ since $ w $ is continuous and positive on $ \overline{\mathcal{C}_{R/2,\tau/2}\times (0,\varrho/2)} $. Thus, we have
\[
L_\alpha (w-\varepsilon \phi)+(w-\varepsilon \phi)_{tt}+a (w-\varepsilon \phi)_t-b (w-\varepsilon \phi)\le0
\]
in $ \mathcal{C}_{R/2,\tau/2}\times (0,\varrho/2) $
with $w-\varepsilon \phi\ge 0 $ on $ \partial (\mathcal{C}_{R/2,\tau/2}\times (0,\varrho/2)) $.
By the weak maximum principle (see Theorem \ref{tw}), we have
\[
w-\varepsilon \phi\ge 0\quad { \rm in}\,\,\mathcal{C}_{R/2,\tau/2}\times (0,\varrho/2).
\]
This implies that
\[
w\ge \varepsilon \phi\ge \varepsilon \theta_{N+1}^{2\alpha}\phi\quad { \rm in}\,\,\mathcal{C}_{R/2,\tau/2}\times (0,\varrho/2).
\]
Consequently, this leads to
\[
\limsup_{\theta_{N+1}\rightarrow 0^+}-\theta_{N+1}^{1-2\alpha}\frac{w(t_0,\tilde{\theta}_0,\theta_{N+1})}{\theta_{N+1}}\le \varepsilon \limsup_{\theta_{N+1}\rightarrow 0^+}-\theta_{N+1}^{1-2\alpha}\frac{\theta_{N+1}^{2\alpha}\phi(t_0)}{\theta_{N+1}}=-\varepsilon\phi(t_0)<0,
\]
as claimed in the proposition. $\Box$

\begin{remark}
 If in addition $ \theta_{N+1}^{1-2\alpha} w_{\theta_{N+1}}\in C(\overline{\mathcal{C}_{R,\tau}\times (0,\varrho)}) $, we have that
 \[
\partial_{\nu^\alpha}w(t_0,\tilde{\theta}_0,0)=- \lim_{\theta_{N+1}\rightarrow 0^+}\theta_{N+1}^{1-2\alpha}\frac{w(t_0, \tilde{\theta}_0,\theta_{N+1})}{\theta_{N+1}}<0.
 \]
\end{remark}

Finally, by the above two Hopf's lemmas and a similar argument as Corollary 4.12 in \cite{CaS}, we obtain the following version strong maximum principle.

\begin{teo}\label{tm}
Let $ \Omega\subset\R\times S^N_+ $ be an open bounded set with a part of boundary $ \Gamma\subset  \R\times \partial S^N_+$. Assume $ w\in C^2(\Omega)\cap C(\bar{\Omega}) $ and $ \theta_{N+1}^{1-2\alpha} w_{\theta_{N+1}}\in C(\bar{\Omega}) $ satisfies
\begin{eqnarray*}
\left\{\begin{array}{l@{\quad }l}
L_\alpha w+w_{tt}+a(t,\theta) w_t-b(t,\theta) w\le0 &{ in}\,\,\Omega,\\
- \lim_{\theta_{N+1}\rightarrow 0^+}\theta_{N+1}^{1-2\alpha}w_{\theta_{N+1}}\ge 0  &{ on}\,\,\Gamma,\\
w\ge0,\quad w\not\equiv 0&{ on}\,\,\Omega,
 \end{array}
 \right.
 \end{eqnarray*}
 where $a$ and $b$ are bounded functions and $b$ is  nonnegative. Then $ w>0 $ in $ \Omega\cup \Gamma $.
\end{teo}

\setcounter{equation}{0}
\section{ Proofs}

We prove our main result in this section via the method of moving planes.  For which we give some preliminary notations, we define
\[\Sigma_\mu=\{(t,\theta):t\in(-\infty,\mu),\theta\in S_+^N\},\]
\[T_\mu=\{(t,\theta):t=\mu,\theta\in S_+^N\},\]
\[w^{\mu}(t,\theta)=\overline{U}(2\mu-t,\theta)-\overline{U}(t,\theta),\]
\[ z^{\mu}(t,\theta)=\overline{V}(2\mu-t,\theta)-\overline{V}(t,\theta).\]
A direct calculation shows the comparison functions $w^{\mu}  $ and $ z^{\mu} $ satisfy
\begin{eqnarray}\label{3.1}
\left\{\begin{array}{l@{\quad }l}
L_\alpha w^{\mu}+w^{\mu}_{tt}+\delta_1 w^{\mu}_t
-\nu_1w^{\mu}=-2\delta_1 \overline{U}_t&{\rm in}\,\,\Sigma_\mu,\\
 -\lim_{\theta_{N+1}\rightarrow 0^+}\theta_{N+1}^{1-2\alpha}\partial_{\theta_{N+1}} w^{\mu}=c^\mu z^{\mu} &{\rm on}\,\,\partial_L \Sigma_\mu,\\
L_\alpha z^{\mu}+z^{\mu}_{tt}+\delta_2 z^{\mu}_t
-\nu_2z^{\mu}=-2\delta_2 \overline{V}_t&{\rm in}\,\,\Sigma_\mu,\\
 -\lim_{\theta_{N+1}\rightarrow 0^+}\theta_{N+1}^{1-2\alpha}\partial_{\theta_{N+1}}  z^{\mu}=d^\mu w^{\mu} &{\rm on}\,\,\partial_L \Sigma_\mu,
 \end{array}
 \right.
 \end{eqnarray}
where $ \partial_L \Sigma_\mu:= (\R\times\partial S^N_+)\cap \overline{\Sigma_\mu} $,
\beq\label{3.2}
c^\mu(t,\theta)=\left\{
  \begin{array}{ll}
\frac{\left(\overline{V}(2\mu-t,\theta)\right)^q-\left(\overline{V}(t,\theta)\right)^q}{\overline{V}(2\mu-t,\theta)-\overline{V}(t,\theta)}   &{\rm if}\,\,\overline{V}(2\mu-t,\theta)\not=\overline{V}(t,\theta),    \\
0   &{\rm if}\,\,\overline{V}(2\mu-t,\theta)=\overline{V}(t,\theta),    \\
\end{array}
\right.
\eqq
and
\beq\label{uu}
d^\mu(t,\theta)=\left\{
  \begin{array}{ll}
\frac{\left(\overline{U}(2\mu-t,\theta)\right)^p-\left(\overline{U}(t,\theta)\right)^p}{\overline{U}(2\mu-t,\theta)-\overline{U}(t,\theta)}   &{\rm if}\,\,\overline{U}(2\mu-t,\theta)\not=\overline{U}(t,\theta),    \\
0   &{\rm if}\,\,\overline{U}(2\mu-t,\theta)=\overline{U}(t,\theta).   \\
\end{array}
\right.
\eqq
From (\ref{3.2}) and (\ref{uu}), we have that
\beq\label{3.3}
c^\mu\ge 0,\quad d^\mu\ge0\quad {\rm in}\,\,\R\times \partial S_+^N.
\eqq
Moreover, the definitions of $ w^\mu $ and $ z^\mu $ imply
\beq\label{3.4}
w^\mu\equiv z^\mu\equiv 0\quad {\rm on}\,\,T_\mu.
\eqq

Next, we show that $ \overline{U} $ and $ \overline{V} $ decay monotonically near $ -\infty $.
In fact, by differentiating (\ref{2.5}), we find that
\beq\label{df}
\overline{U}_t=r^{\beta_1}(\beta_1 U+rU_r)\quad{\rm and}\quad\overline{V}_t=r^{\beta_2}(\beta_2 V+rV_r).
\eqq
So take into account $ \beta_1,\beta_2>0 $, $ r=e^t $, $ \overline{U}>0 $ and $ \overline{V}>0 $, we can obtain $ t_0 $ for which
\beq\label{3.7}
\overline{U}_t>0\quad {\rm and}\quad \overline{V}_t>0\quad {\rm in }\,\,\Sigma_{t_0},
\eqq
and
\beq\label{3.8}
0<\overline{U}(t,\theta), \overline{V}(t,\theta)<\varepsilon_0\quad {\rm in }\,\,\Sigma_{t_0},
\eqq
where $ 0<\varepsilon_0<<1$.

Now, we fix $ t_0 $ such that (\ref{3.7}) and (\ref{3.8}) holds. The following maximum principle for system (\ref{3.1}) near $ -\infty $ is needed, which permits us to get the moving planes method started.

\begin{lemma}\label{l1}

(1) For all $ \mu\in(-\infty, t_0] $ one has
$ w^\mu\ge0 $ and $ z^\mu\ge0 $ in $ \Sigma_{\mu}$. \\
(2)  Suppose that for $ \mu\in( t_0,+\infty) $, we have
$ w^\mu\ge0 $ and $ z^\mu\ge0 $ on  $ T_{t_0}$.
Then $ w^\mu\ge0 $ and $ z^\mu\ge0 $ in  $\Sigma_{t_0}$.
\end{lemma}

{\bf Proof.}
Observe that in both cases $ (1) $ and $ (2) $ we have  $w^\mu\ge0 $ and $ z^\mu\ge0 $ on  $T_{t_0\wedge \mu}$ by (\ref{3.4}),  where $ t_0\wedge \mu=\min\{t_0,\mu\} $.  Therefore, we treat both cases at the same time by a contradiction argument,
 assuming that
\beq\label{3.15}
\min\left\{\inf_{\Sigma_{ {t_0\wedge \mu}}}w^\mu,\,\, \inf_{\Sigma_{ {t_0\wedge \mu}}}z^\mu\right\}<0.
\eqq

Up to an symmetry in the argument, there are two cases to rule out.

Case I:
\[\inf_{\Sigma_{ {t_0\wedge \mu}}}w^\mu<0\quad{\rm and}\quad \inf_{\Sigma_{{t_0\wedge \mu}}}z^\mu\ge0.\]

We consider the function
\begin{equation}\label{3.33}
W(x) = \left\{ \arraycolsep=1.5pt
\begin{array}{lll}
\max\{-w^\mu(x),0\}, & \quad  x\in\overline{\Sigma_{ {t_0\wedge \mu}}},\\[2mm]
0,   & \quad  x\in\left(\overline{\Sigma_{ {t_0\wedge \mu}}}\right)^c.
 \end{array}
 \right.
 \end{equation}
Therefore, using the equation (\ref{3.1}), we have
\begin{eqnarray}\nonumber
0&\le&\int_{\Sigma_{ {t_0\wedge \mu}}} \theta_{N+1}^{1-2\alpha}|\nabla W|^2 e^{\delta_1t}d\theta dt=-\int_{\Sigma_{ {t_0\wedge \mu}}} \theta_{N+1}^{1-2\alpha}\nabla w^\mu\nabla W\cdot e^{\delta_1t}d\theta dt\\\nonumber
&=&\int_{\Sigma_{ {t_0\wedge \mu}}} div(\theta_{N+1}^{1-2\alpha}\nabla w^\mu) We^{\delta_1t}d\theta dt-\int_{\partial_L\Sigma_{ {t_0\wedge \mu}}} \theta_{N+1}^{1-2\alpha} \frac{\partial w^\mu}{\partial \nu} We^{\delta_1t}d\theta dt\\\nonumber
&=&-\int_{\Sigma_{ {t_0\wedge \mu}}} \theta_{N+1}^{1-2\alpha}(w^\mu_{tt}+\delta_1 w^\mu_t) We^{\delta_1t}d\theta dt+\int_{\Sigma_{ {t_0\wedge \mu}}} \theta_{N+1}^{1-2\alpha}\nu_1 w^\mu We^{\delta_1t}d\theta dt\\\nonumber
&-&\int_{\Sigma_{ {t_0\wedge \mu}}} \theta_{N+1}^{1-2\alpha} 2\delta_1\overline{U}_tWe^{\delta_1t}d\theta dt-\int_{\partial_L\Sigma_{ {t_0\wedge \mu}}} c^\mu z^\mu We^{\delta_1t}d\theta dt\\\label{3.37}
&\le&-\int_{\Sigma_{ {t_0\wedge \mu}}} \theta_{N+1}^{1-2\alpha}(w^\mu_{tt}+\delta_1 w^\mu_t) We^{\delta_1t}d\theta dt
\end{eqnarray}
since (\ref{3.7}), $\delta_1\ge0$, $\nu_1>0$ and $z^\mu\ge0$ on $\partial_L\Sigma_{ {t_0\wedge \mu}}$ by the continuity of $z^\mu$.

Since $ \beta_1>0 $, the definitions of $ \overline{U} $ (see (\ref{2.5})) and $ \overline{w^\mu} $ imply
\[
\liminf_{t\rightarrow-\infty} \inf_{\theta\in S_+^N} w^\mu(t,\theta)\ge0.
\]
This implies that $W\equiv0$ as $t\rightarrow-\infty$ and $\theta\in S_+^N$. Therefore, we can estimate that
\begin{eqnarray}\nonumber
\int_{\Sigma_{ {t_0\wedge \mu}}} \theta_{N+1}^{1-2\alpha}(w^\mu_{tt}+\delta_1 w^\mu_t) We^{\delta_1t}d\theta dt&=& \int_{\Sigma_{ {t_0\wedge \mu}}} \theta_{N+1}^{1-2\alpha}\left(e^{\delta_1t} w^\mu_t\right)_t Wd\theta dt\\\nonumber
&=& -\int_{\Sigma_{ {t_0\wedge \mu}}} \theta_{N+1}^{1-2\alpha}e^{\delta_1t} w^\mu_t W_td\theta dt\\\nonumber
&=& \int_{\Sigma_{ {t_0\wedge \mu}}} \theta_{N+1}^{1-2\alpha}e^{\delta_1t} |W_t|^2d\theta dt\\\label{3.38}
&\ge&0.
\end{eqnarray}
Together (\ref{3.37}) and (\ref{3.38}), we have
\[
\int_{\Sigma_{ {t_0\wedge \mu}}} \theta_{N+1}^{1-2\alpha}|\nabla W|^2 e^{\delta_1t}d\theta dt\equiv 0.
\]
This is impossible unless $W\equiv 0$ in $\Sigma_{ {t_0\wedge \mu}}$ and therefore $w^\mu\ge0$ in $\Sigma_{ {t_0\wedge \mu}}$. This contradicts with $\inf_{\Sigma_{ {t_0\wedge \mu}}}w^\mu<0$.

 Case II:
\[\inf_{\Sigma_{ {t_0\wedge \mu}}}w^\mu<0\quad{\rm and}\quad \inf_{\Sigma_{ {t_0\wedge \mu}}}z^\mu<0.\]

Since $ \overline{w^\mu}\ge0 $ on $ T_{t_0\wedge \mu} $ and $\liminf_{t\rightarrow-\infty} \inf_{\theta\in S_+^N} w^\mu(t,\theta)\ge0$, there exists a point
$ (\bar{t},\bar{\theta})\in \Sigma_{ {t_0\wedge \mu}}\cup \partial_L  \Sigma_{ {t_0\wedge \mu}}$
such that the negative infimum of $ w^\mu $ is achieved, that is,
$$ w^\mu(\bar{t},\bar{\theta})=\inf_{\Sigma_{ {t_0\wedge \mu}}}w^\mu<0.$$
If $ (\bar{t},\bar{\theta})\in\Sigma_{ {t_0\wedge \mu}} $, then
\[w^\mu_{tt}(\bar{t},\bar{\theta})\ge0,\quad \Delta_\theta w^\mu(\bar{t},\bar{\theta})\ge0,\quad  w^\mu_t(\bar{t},\bar{\theta})=0,\quad{\rm and}\quad\nabla w^\mu(\bar{t},\bar{\theta})=0. \]
Thus, from (\ref{2.8}), (\ref{3.7}) and the first equation of (\ref{3.1}), we have that
\[0<-\nu_1w^\mu(\bar{t},\bar{\theta})\le-2\delta_1 \overline{U}_t(\bar{t},\bar{\theta})\le 0,\]
since $ \nu_1>0 $ and this is a contradiction.
So $ (\bar{t},\bar{\theta})\in\partial_L \Sigma_{t_0\wedge \mu}$, which implies that $  \partial_{ \theta_{N+1}} \overline{w^\mu}(\bar{t},\bar{\theta})\ge0 $. Therefore, $c^\mu(\bar{t},\bar{\theta})z^{\mu}(\bar{t},\bar{\theta})\le0$ and thus $z^{\mu}(\bar{t},\bar{\theta})\le0$ since $ c^\mu\ge0 $. By the continuity of $z^\mu$, we know $z^{\mu}(t,\theta)\le0$ on
$\partial_L \Sigma_{t_0\wedge \mu}$.
Then
\[
0<V(2\mu-t,\theta)\le V(t,\theta)<\varepsilon_0\quad {\rm on }\,\,\partial_L \Sigma_{t_0\wedge \mu}.
\]
By the mean value principle, we have $c^\mu(t,\theta)\le q \varepsilon_0^{q-1}$ on $\partial_L \Sigma_{t_0\wedge \mu}$. Similarly, we have $d^\mu(t,\theta)\le p\varepsilon_0^{p-1}$ on $\partial_L \Sigma_{t_0\wedge \mu}$.

We define function
\begin{equation}\label{3.34}
Z(x) = \left\{ \arraycolsep=1.5pt
\begin{array}{lll}
\max\{-z^\mu(x),0\}, & \quad  x\in\overline{\Sigma_{ {t_0\wedge \mu}}},\\[2mm]
0,   & \quad  x\in\left(\overline{\Sigma_{ {t_0\wedge \mu}}}\right)^c.
 \end{array}
 \right.
 \end{equation}

Without of loss generality, we suppose $\delta_1\le \delta_2$. By a similar estimate in Case I, we have
\begin{eqnarray*}
\int_{\Sigma_{ {t_0\wedge \mu}}}  \theta_{N+1}^{1-2\alpha}|\nabla Z|^2 e^{\delta_1t}d\theta dt
&\le&\int_{\Sigma_{ {t_0\wedge \mu}}}  \theta_{N+1}^{1-2\alpha}|\nabla Z|^2 e^{\delta_2t}d\theta dt\\
&\le& -\int_{\partial_L\Sigma_{ {t_0\wedge \mu}}} d^\mu w^\mu Ze^{\delta_2t}d\theta dt\\
&=&\int_{\partial_L\Sigma_{ {t_0\wedge \mu}}} d^\mu W Ze^{\delta_2t}d\theta dt\\
&=&\int_{\partial_L\Sigma_{ {t_0\wedge \mu}}} d^\mu W Ze^{(\delta_2-\delta_1)t}e^{\delta_1t}d\theta dt\\
&\le& p\varepsilon_0^{p-1}e^{(\delta_2-\delta_1)t_0}\int_{\partial_L\Sigma_{ {t_0\wedge \mu}}}  W Ze^{\delta_1t}d\theta dt.
\end{eqnarray*}
By the H\"older and fractional Sobolev trace inequalities as in \cite{T} (see also \cite{X} for the Sobolev  trace inequality in all $\R^n$), we know
\begin{flushleft}
$\int_{\partial_L\Sigma_{ {t_0\wedge \mu}}}  W Ze^{\delta_1t}d\tilde{\theta} dt$\\
$\le\left(\int_{\partial_L\Sigma_{ {t_0\wedge \mu}}}  |e^{\delta_1t/2}Z|^2d\tilde{\theta} dt\right)^{1/2} \left(\int_{\partial_L\Sigma_{ {t_0\wedge \mu}}}  |e^{\delta_1t/2}W|^2d\tilde{\theta} dt\right)^{1/2}$\\
$\le C_{N,\alpha}\left(\int_{\Sigma_{ {t_0\wedge \mu}}}  \theta_{N+1}^{1-2\alpha}|\nabla Z|^2 e^{\delta_1t}d\theta dt\right)^{1/2} \left(\int_{\Sigma_{ {t_0\wedge \mu}}}  \theta_{N+1}^{1-2\alpha}|\nabla W|^2 e^{\delta_1t}d\theta dt\right)^{1/2},$
\end{flushleft}
where $\theta=(\tilde{\theta},\theta_{N+1})\in S^N_+$ and $C_{N,\alpha}$ is a positive constant depending only on $N$ and $\alpha$.
Hence,
\beq\label{3.35}
\left(\int_{\Sigma_{ {t_0\wedge \mu}}}  \theta_{N+1}^{1-2\alpha}|\nabla Z|^2 e^{\delta_1t}d\theta dt\right)^{1/2}\le p\varepsilon_0^{p-1} C_{N,\alpha}e^{(\delta_2-\delta_1)t_0}\left(\int_{\Sigma_{ {t_0\wedge \mu}}}  \theta_{N+1}^{1-2\alpha}|\nabla W|^2 e^{\delta_1t}d\theta dt\right)^{1/2}.
\eqq
Moreover, as the argument in Case I, we have
\[
\int_{\Sigma_{ {t_0\wedge \mu}}} \theta_{N+1}^{1-2\alpha}|\nabla W|^2 e^{\delta_1t}d\theta dt\le -\int_{\partial_L\Sigma_{ {t_0\wedge \mu}}} c^\mu z^\mu We^{\delta_1t}d\tilde{\theta} dt
=\int_{\partial_L\Sigma_{ {t_0\wedge \mu}}} c^\mu Z We^{\delta_1t}d\tilde{\theta} dt.
\]
Similarly, by the H\"older and fractional Sobolev trace inequalities, we have
\begin{flushleft}
 $\int_{\Sigma_{ {t_0\wedge \mu}}} \theta_{N+1}^{1-2\alpha}|\nabla W|^2 e^{\delta_1t}d\theta dt$\\
 $\le q\varepsilon_0^{q-1}\left(\int_{\partial_L\Sigma_{ {t_0\wedge \mu}}}  |e^{\delta_1t/2}Z|^2d\tilde{\theta} dt\right)^{1/2} \left(\int_{\partial_L\Sigma_{ {t_0\wedge \mu}}}  |e^{\delta_1t/2}W|^2d\tilde{\theta} dt\right)^{1/2}$\\
$\le q\varepsilon_0^{q-1}C_{N,\alpha}\left(\int_{\Sigma_{ {t_0\wedge \mu}}}  \theta_{N+1}^{1-2\alpha}|\nabla Z|^2 e^{\delta_1t}d\theta dt\right)^{1/2} \left(\int_{\Sigma_{ {t_0\wedge \mu}}}  \theta_{N+1}^{1-2\alpha}|\nabla W|^2 e^{\delta_1t}d\theta dt\right)^{1/2}.$
\end{flushleft}
Therefore,
\beq\label{3.36}
\left(\int_{\Sigma_{ {t_0\wedge \mu}}}  \theta_{N+1}^{1-2\alpha}|\nabla W|^2 e^{\delta_1t}\right)^{1/2}d\theta dt\le q\varepsilon_0^{q-1} C_{N,\alpha}\left(\int_{\Sigma_{ {t_0\wedge \mu}}}  \theta_{N+1}^{1-2\alpha}|\nabla Z|^2 e^{\delta_1t}d\theta dt\right)^{1/2}.
\eqq
Then, combining (\ref{3.35}) and (\ref{3.36}), we have
\begin{flushleft}
$\left(\int_{\Sigma_{ {t_0\wedge \mu}}}  \theta_{N+1}^{1-2\alpha}|\nabla W|^2 e^{\delta_1t}d\theta dt\right)^{1/2}$\\
$\le pq C_{N,\alpha}e^{(\delta_2-\delta_1)t_0}\varepsilon_0^{p+q-2}\left(\int_{\Sigma_{ {t_0\wedge \mu}}}  \theta_{N+1}^{1-2\alpha}|\nabla W|^2e^{\delta_1t}d\theta dt\right)^{1/2}$
\end{flushleft}
and
\begin{flushleft}
$\left(\int_{\Sigma_{ {t_0\wedge \mu}}}  \theta_{N+1}^{1-2\alpha}|\nabla Z|^2 e^{\delta_1t}d\theta dt\right)^{1/2}$\\
$\le pq C_{N,\alpha}e^{(\delta_2-\delta_1)t_0}\varepsilon_0^{p+q-2}\left(\int_{\Sigma_{ {t_0\wedge \mu}}}  \theta_{N+1}^{1-2\alpha}|\nabla Z|^2e^{\delta_1t}d\theta dt\right)^{1/2}.$
\end{flushleft}
These are impossible since $\varepsilon_0<<1$ and $p+q>2$ since $pq>1$ unless $W\equiv0$ and $Z\equiv0$ in $\Sigma_{ {t_0\wedge \mu}}$ and thus $w^\mu\ge0$ and $z^\mu\ge0$ in $\Sigma_{ {t_0\wedge \mu}}$, which contradict with our assumption.
 We complete the proof of Lemma \ref{l1}. $ \Box $\\

By Case (1) of Lemma \ref{l1}, we have $ w^\mu\ge0 $ and $ z^\mu\ge0 $ in $ \Sigma_{\mu} $ for all $ \mu\in(-\infty,t_0] $. This enables us to define a maximal value of $ \mu $ up to which the positivity of these functions holds. This is the purpose of the following lemma.

\begin{lemma}\label{l2}
We have either
\begin{flushleft}\label{3.21}
$ (1)$: there exists $\bar{\mu}\in\R$ such that $w^{\bar{\mu}}\equiv0$ and $z^{\bar{\mu}}\equiv0$ in $\Sigma_{\bar{\mu}}$, or
\end{flushleft}

\begin{flushleft}\label{3.22}
$(2)$: for any $\mu\in\R$ one has $w^{\mu}>0$ and $z^{\mu}>0$ in $\Sigma_{\mu}$.
\end{flushleft}
Moreover, in the latter case one has
\beq\label{3.23}
\overline{U}_t>0\quad and \quad \overline{V}_t>0 \quad in\,\, \R\times S_+^N.
\eqq
\end{lemma}

{\bf Proof.}
Define
\beq\label{3.24}
\Lambda=\sup\{\mu\in\R: \forall \lambda\in(-\infty,\mu),\,\, w^{\lambda}\ge0\,\, {\rm and}\,\,z^{\lambda}\ge0\,\, {\rm in}\,\, \Sigma_{\lambda}\}.
\eqq
By Lemma \ref{l1}, it is clear that $ \Lambda>-\infty $.
Next, we prove that either $ \Lambda<+\infty $, in which case (1) holds with $ \bar{\mu}=\Lambda $, or $ \Lambda=+\infty $ in which case  we have case (2) together with
(\ref{3.23}).

If $\Lambda=+\infty $, then case (2) is trivially satisfied. Moreover, in this case (\ref{3.23}) is a consequence of the following argument. Since (\ref{3.4}), we have
\[
\frac{\partial w^\mu}{\partial t}\le0\quad{\rm and}\quad \frac{\partial z^\mu}{\partial t}\le0\quad{\rm in}\,\,\Sigma_{\mu}.
\]
 By the definitions of $w^\mu$ and $z^\mu$, we know
 \beq\label{3.40}
 \frac{\partial w^\mu}{\partial t}=-2\overline{U}_t\quad{\rm and}\quad \frac{\partial z^\mu}{\partial t}-2\overline{V}_t\quad{\rm on}\,\,T_{\mu}.
 \eqq
 Therefore,
 \[
 \overline{U}_t\ge0 \quad{\rm and}\quad \overline{V}_t\ge0\quad{\rm on}\,\,T_{\mu}
 \]
 for all $\mu\in\R$ and thus throughout $\R\times S_+^N$.
 Then, by (\ref{3.1}), we have
 \begin{eqnarray}\label{3.39}
\left\{\begin{array}{l@{\quad }l}
L_\alpha w^{\mu}+w^{\mu}_{tt}+\delta_1 w^{\mu}_t
-\nu_1w^{\mu}\le0&{\rm in}\,\,\Sigma_\mu,\\
L_\alpha z^{\mu}+z^{\mu}_{tt}+\delta_2 z^{\mu}_t
-\nu_2z^{\mu}\le0&{\rm in}\,\,\Sigma_\mu.
 \end{array}
 \right.
 \end{eqnarray}
Applying the Hopf lemma to each equation in (\ref{3.39}) yields
\[
\frac{\partial w^\mu}{\partial t}<0\quad{\rm and}\quad \frac{\partial z^\mu}{\partial t}<0\quad{\rm on}\,\,T_{\mu},
\]
and thus we have (\ref{3.23}) thanks to (\ref{3.40}).

 Suppose that $ \Lambda<+\infty $. We prove case (1) by contradiction and assume that $w^{\Lambda}\not\equiv0$  or $z^{\Lambda}\not\equiv0$ in\,\, $\Sigma_{\Lambda}$. For all $ -\infty<\mu\le\Lambda $, by the stronger maximum principle, we know that
$w^{\mu}>0$  and  $z^{\mu}>0$ in $\Sigma_{\mu}$. The above arguments imply that $ \overline{U}_t>0 $ and $ \overline{V}_t>0 $
on  $T_{\mu}$ for $ -\infty<\mu<\Lambda $.  Hence, $ \overline{U}_t>0 $ and $ \overline{V}_t>0 $ in  $\Sigma_{\Lambda}$.

Therefore, by (\ref{3.1}) we have that, for $ -\infty<\mu<\Lambda $,
\begin{eqnarray}\label{3.25}
\left\{\begin{array}{l@{\quad }l}
L_\alpha w^{\mu}+w^{\mu}_{tt}+\delta_1 w^{\mu}_t
-\nu_1w^{\mu}=-2\delta_1 \overline{U}_t\le0&{\rm in}\,\,\Sigma_\mu,\\
L_\alpha z^{\mu}+z^{\mu}_{tt}+\delta_2 z^{\mu}_t
-\nu_2z^{\mu}=-2\delta_2 \overline{V}_t\le0&{\rm in}\,\,\Sigma_\mu.\\
 \end{array}
 \right.
 \end{eqnarray}
Now, evaluating (\ref{3.25}) at $ \mu=\Lambda $ by continuity, we obtain
\beq\label{3.26}
L_\alpha w^{\Lambda}+w^{\Lambda}_{tt}+\delta_1 w^{\Lambda}_t
-\nu_1w^{\Lambda}\le0\quad{\rm in}\,\,\Sigma_\Lambda,
 \eqq
and
\beq\label{3.27}
L_\alpha z^{\Lambda}+z^{\Lambda}_{tt}+\delta_2 z^{\Lambda}_t
-\nu_2z^{\Lambda}\le0\quad{\rm in}\,\,\Sigma_\Lambda.
 \eqq
An application of the strong maximum principle (see Theorem \ref{tm}) to (\ref{3.26}) and (\ref{3.27}) implies that either
$  w^{\Lambda}>0 $  or $  w^{\Lambda}\equiv0 $ in $ \Sigma_\Lambda$ on the one hand $  z^{\Lambda}>0 $  or $  z^{\Lambda}\equiv0 $ in $ \Sigma_\Lambda$ on the other hand.

It is easy to check that (\ref{3.1}) and Theorem \ref{tm} rules out the cases $  w^{\Lambda}>0 $ and  $  z^{\Lambda}\equiv0 $ in $ \Sigma_\Lambda$ as well as  $  w^{\Lambda}\equiv0 $ and $  z^{\Lambda}>0 $
 in $ \Sigma_\Lambda$. Here we only show the case $  w^{\Lambda}>0 $ and  $  z^{\Lambda}\equiv0 $ in $ \Sigma_\Lambda$ is impossible. In fact, since  $  z^{\Lambda}\equiv0 $ in $ \Sigma_\Lambda$, by the continuous
 up to the boundary, we have $  z^{\Lambda}=0 $  and $\theta_{N+1}^{1-2\alpha}\partial_{\theta_{N+1}}  z^{\Lambda}=0$ on $ \partial_L \Sigma_\Lambda$. If
$  w^{\Lambda}>0 $ in $ \Sigma_\Lambda$, then applying Theorem \ref{tm} to (\ref{3.1}) we  know $  w^{\Lambda}>0 $ on $ \partial_L\Sigma_\Lambda$ and thus $ d^\Lambda>0 $. So by (\ref{3.1}) we have
 \beq\label{44}
 0=-\lim_{\theta_{N+1}\rightarrow 0^+}\theta_{N+1}^{1-2\alpha}\partial_{\theta_{N+1}}  z^{\Lambda}=d^\Lambda(t,\theta)w^{\Lambda}>0
 \eqq
 on $\partial_L\Sigma_\Lambda$. This is a contradiction.

Hence in the following we may assume that both $  w^{\Lambda} $ and $  z^{\Lambda} $
are strictly positive in $ \Sigma_\Lambda$. By the Hopf's lemma (see Lemma \ref{h1}) we have
\beq\label{h12}
\frac{\partial w^{\Lambda}}{\partial t}<0 \quad {\rm and} \quad \frac{\partial z^{\Lambda}}{\partial t}<0 \quad {\rm on}\,\, T_\Lambda,
\eqq
since $ w^{\Lambda}=0$ and $ z^{\Lambda}=0$ on $T_\Lambda$.

Next, we claim that there exists $ \varepsilon>0 $ such that $  w^{\mu}\ge0 $ and $  z^{\mu}\ge0 $
in $ \Sigma_\mu$ for all $ \mu\in(\Lambda,\Lambda+\varepsilon) $. This will result in a contradiction with the definition of $ \Lambda $; hence the lemma will be proved. This is done in the following way.

We split the domain into three disjoint subsets:
\[\Sigma_\mu=\Sigma_{t_0}\cup(\overline{\Sigma_{\Lambda-\delta}}\setminus \Sigma_{t_0})\cup (\Sigma_\mu\setminus\overline{\Sigma_{\Lambda-\delta}}),\]
for some small $ \delta >0$ to be defined.

We start by checking the second set. For a given $ \delta >0$, we know $  w^{\Lambda}>0 $ and $  z^{\Lambda}>0 $ in the compact set $ \overline{\Sigma_{\Lambda-\delta}}\setminus \Sigma_{t_0} $. Therefore, a straightforward continuity argument implies that
there exists $ \varepsilon_2=\varepsilon_2(\delta)>0$ such that  for all $ \mu\in[\Lambda,\Lambda+\varepsilon_2) $,
\[
\min\left\{\inf_{\overline{\Sigma_{\Lambda-\delta}}\setminus \Sigma_{t_0}}w^\mu,\inf_{\overline{\Sigma_{\Lambda-\delta}}\setminus \Sigma_{t_0}}z^\mu\right\}>0.
\]

We carry on the analysis by examining the first part of the domain. By the above consideration, for all $ \mu\in[\Lambda,\Lambda+\varepsilon_2) $, we have $  w^{\mu}\ge0 $ and $  z^{\mu}\ge0 $ on $ T_{t_0}\subset (\overline{\Sigma_{\Lambda-\delta}}\setminus \Sigma_{t_0})$. An application of  case (2) of Lemma \ref{l1}, we have $  w^{\mu}>0$ and $  z^{\mu}>0$ in $ \Sigma_{t_0} $.

Finally, we do the analysis on the third  part of the domain, namely $ \Sigma_\mu\setminus \overline{\Sigma_{\Lambda-\delta}} $.
A simple continuity argument shows that (\ref{h12}) remains valid if $ \Lambda $ is replaced by any $ \mu $ in a small right neighborhood of $ \Lambda $, that is, there exists $ \varepsilon_3>0 $ such that for all $ \mu\in[\Lambda,\Lambda+\varepsilon_3) $,
\[
\frac{\partial w^\mu}{\partial t}<0\quad{\rm and}\quad  \frac{\partial z^\mu}{\partial t}<0\quad {\rm on}\,\,T_\mu.
\]
By elliptic estimates give locally uniform $C^{2}$ bounds for $w^\mu$  and  $z^\mu$ in $(t,\theta)$ as well as in $\mu$. Hence
\[
\inf_{\mu-\varepsilon<t<\mu}\left(-\frac{\partial w^\mu}{\partial t}\right)\ge \inf_{T_\mu}\left(-\frac{\partial w^\mu}{\partial t}\right)-C\varepsilon,
\]
for any $\varepsilon\in(0,1)$ and $\theta\in S^N_+$. Similarly for $z^\mu$. This means that we can choose $\varepsilon_3^\prime\in (0,\varepsilon_3)$ such that $\partial w^\mu/\partial t<0$ and $\partial z^\mu/\partial t<0$ in the set $ \{(t,\theta):\mu-\varepsilon_3^\prime<t<\mu,\theta\in S^N_+ \} $ for all $\mu\in [\Lambda,\Lambda+\varepsilon_3^\prime/2)$.
Hence, by $  w^{\mu}=0 $ and $  z^{\mu}=0 $ on $ T_\mu$, we know
 $ w^\mu>0 $ and $ z^\mu>0 $ in the set $ \{(t,\theta):\mu-\varepsilon_3^\prime<t<\mu,\theta\in S^N_+ \} $ for all $\lambda\in [\Lambda,\Lambda+\varepsilon_3^\prime/2)$.

Fixing now $ \delta=\varepsilon_3^\prime/2 $ and $ \varepsilon=\min\{\varepsilon_2(\delta),\varepsilon_3\} $, summing up the above results,
we have $ w^\mu >0$ and $ z^\mu >0$ in $ \Sigma_\mu$ for all $ \mu\in(\Lambda,\Lambda+\varepsilon) $. This contradicts to the definition of $ \Lambda $ (see (\ref{3.24})). $ \Box $\\

{\bf Proof of Theorem \ref{t1}.} Suppose that $ (U,V) $ is positive solution to (\ref{ex}).
Then the comparison functions $ w^\mu $ and $ z^\mu $  satisfy the alternative in Lemma \ref{l1}. Next, we prove both cases $(1)$ and $(2)$ in Lemma \ref{l2} cannot happen.

First, we show the case $(1)$ in Lemma \ref{l2} provides a contradiction. In order to get the contradiction, we translate the origin to $ \bar{\mu} $, that is,
define $ \hat{U}(t,\theta)=\overline{U}(t+\bar{\mu},\theta) $ and $ \hat{V}(t,\theta)=\overline{V}(t+\bar{\mu},\theta) $.
Since $w^{\bar{\mu}}\equiv0$ and $z^{\bar{\mu}}\equiv0$ in $\Sigma_{\bar{\mu}}$, then those two functions are even in the variable $ t $, i.e.,
\beq\label{3.30}
\hat{U}(-t,\theta)=\hat{U}(t,\theta)\quad{\rm and}\quad \hat{V}(-t,\theta)=\hat{V}(t,\theta),
\eqq
This implies $\overline{U}_t  $ and $\overline{V}_t  $ are odd functions. On the other hand, by the first and  third equations of (\ref{3.1})
and $ (\delta_1,\delta_2)\neq(0,0) $ (see \ref{2.8}), we know $\overline{U}_t  $ or $\overline{V}_t  $ are even in variable $ t $. So we can conclude that
$ \overline{U} $ or $ \overline{V} $ must be constant in the whole domain $ \R\times S_+^N $. This contradicts the regularity of $ \overline{U} $ or $ \overline{V} $  at the origin, see (\ref{2.5}).

We complete the proof of Theorem \ref{t1} by showing that case $(2)$ in Lemma \ref{l2} is also impossible. Observe that (\ref{ex}) is translation invariant in $ x $ direction for $ X=(x,y) \in \R^N\times\R^+$. Hence we can change the initial function, i.e., we define $ U^{x_0}(x,y)=U(x-x_0,y) $, $ V^{x_0}(x,y)=V(x-x_0,y) $ (and corresponding $ \overline{U}^{x_0} $, $ \overline{V}^{x_0} $) for any $ x_0\in\R^N $. Repeat the whole discussion for these new functions, only two cases can be arise. First, there exists an origin $ x_0 $ such that case $(1)$ holds. As we have done, it is impossible. Another case is, for any origin $ x_0\in\R^N $, $(2)$ holds for the transformed  functions $ \overline{U}^{x_0} $ and $ \overline{V}^{x_0} $. So by (\ref{df}), we have in particular
\[\beta_1 U(X)+ \nabla U(X)\cdot (x-x_0,y)\ge0\]
 for any $X=(x,y)\in\R^N\times \R^+$ and $x_0\in\R^N$.
Hence we can obtain
\beq\label{3.31}
\nabla U(X)\cdot (\frac{x-x_0}{|(x-x_0,y)|},\frac{y}{|(x-x_0,y)|})\ge -\frac{\beta_1 U(X)}{|(x-x_0,y)|}
\eqq
for all $ X\in (\R^{N}\setminus \{x_0\})\times\R^+ $. Now we let $ e=(e_1,e_2,\cdots,e_{N+1})\in S_+^N $ with $ e_{N+1}>0 $ fixed and choose $x_0=x-\sigma (e_1,\cdots,e_{N})$ and $ y=\sigma e_{N+1} $
 for $ \sigma>0 $. Thus we rewrite (\ref{3.31}) as
 \[\nabla U(X)\cdot e\ge -\frac{\beta_1 U(X)}{\sigma}.\]
Letting  $ \sigma\rightarrow+\infty $ yields
\beq\label{3.32}
\nabla U(X)\cdot e\ge 0.
\eqq
Since  (\ref{3.32}) holds for any $ e\in S_+^N $ and $ X\in\R^{N+1}_+ $, then we deduce that
\[
\nabla U(X)\equiv 0.
\]
This is impossible by the second equation of (\ref{ex}) since the solution we are dealing with is nontrivial. $ \Box $

\setcounter{equation}{0}
\section{ Acknowledgements}
A. Q. was partially supported by Fondecyt Grant
No. 1151180 Programa Basal, CMM. U. de Chile and Millennium Nucleus
Center for Analysis of PDE NC130017 and
A. Xia was supported by the "Programa de Iniciaci\'on a la  Investigaci\'on Cient\'ifica" (PIIC) UTFSM 2014.

\end{document}